\documentclass[a4paper,12pt]{amsart}
\usepackage{datetime}

\usepackage{amssymb}
\usepackage{amsfonts}

\newcommand{\parq}{\partial_q}

 \usepackage{amsmath}
\usepackage{amsthm}
\newtheorem{thm}{Theorem}[section]

\newtheorem{lem}{Lemma}[section]

\usepackage{graphicx}
\usepackage[T1]{fontenc} %skandit
\usepackage{tabularx}
\usepackage{multirow}
\usepackage{color,soul,enumerate}
\usepackage{tikz}
\usetikzlibrary{arrows.meta, calc, decorations.markings}
\usepackage{amsfonts}
\usepackage{pgfplots}
\usepackage{float}

\theoremstyle{definition}

\newtheorem{example}{Example}[section]

\newtheorem{rem}{Remark}[section]

\usepackage{hyperref}
\hypersetup{
    colorlinks=true,
    citecolor=blue,
    linkcolor=blue,
    urlcolor=blue
}
\nonstopmode \numberwithin{equation}{section}
\newcounter {own}
\def\theown {\thesection       .\arabic{own}}

{\qed\bigskip}

\newcounter{alphabet}
\newcounter{tmp}

\newtheorem{innercustomthm}{Theorem}
\newenvironment{customthm}[1]
  {\renewcommand\theinnercustomthm{#1}\innercustomthm}
  {\endinnercustomthm}
\makeatletter
\def\thefootnote{\@arabic\c@footnote}
\makeatother
\renewcommand{\thefootnote}{\arabic{footnote}}

\begin{document}
\title{The $q$-Jack Lemma and Its Applications}

\author{Md Firoz Ali}
\address{Md Firoz Ali,
 National Institute Of Technology Durgapur, West Bengal, India}
\email{fali.maths@nitdgp.ac.in, ali.firoz89@gmail.com}

\author{Lokenath Thakur}
\address{Lokenath Thakur, National Institute Of Technology Durgapur, West Bengal, India}
\email{lokenaththakur1729@gmail.com}


%\def\thefootnote{}
%\footnotetext{ {\tiny File:~\jobname.tex,
%printed: \number\year-\number\month-\number\day,
 %         \thehours.\ifnum\theminutes<10{0}\fi\theminutes }
%} \makeatletter\def\thefootnote{\@arabic\c@footnote}\makeatother

\begin{abstract}
In this paper, we establish a Jackson $q$-difference analogue of Jack's lemma and its higher-order extension by Miller and Mocanu, which serve as cornerstone tools in geometric function theory and differential subordination. Applying this foundation, we establish $q$-counterparts of Nunokawa's boundary theorem and the Miller-Mocanu subordination lemma, and deduce applications to circular enclosure bounds.  All obtained results naturally recover their classical counterparts in the limiting case as $q \to 1^-$.
\end{abstract}
\maketitle
\thanks{}

\maketitle
\pagestyle{myheadings}
\markboth{Md Firoz Ali and Lokenath Thakur}{The $q$-Jack Lemma and Its Applications}

\section{Introduction}
Let $\mathcal{H}$ denote the class of all analytic functions in the unit disk $\mathbb{D}=\{z\in\mathbb{C}:\vert{}z\vert{}<1\}$, and for any $r > 0$, let $\mathbb{D}_r = \{z \in \mathbb{C} : \vert{}z\vert{} < r\}$. Let $\mathcal{A}$ be the subclass of $\mathcal{H}$ consisting of functions $f$ of the form
\begin{align}\label{FL-001}
f(z)=z+\sum_{n=2}^\infty a_n z^n.
\end{align}
Further, let $\mathcal{S}$ be the subclass of $\mathcal{A}$ consisting of functions that are univalent in $\mathbb{D}$. A function $f \in \mathcal{A}$ is said to be starlike or convex if $f$ maps $\mathbb{D}$ conformally onto a domain that is starlike with respect to the origin or convex, respectively. The subclasses of $\mathcal{S}$ consisting of starlike and convex functions are denoted by $\mathcal{S}^*$ and $\mathcal{K}$, respectively. It is well-established that $\mathcal{K} \subset \mathcal{S}^* \subset \mathcal{S}$.\\

The study of $q$-calculus, also called quantum calculus, goes back to Euler, who in the eighteenth century used $q$-series in his work on partitions, and to Gauss and Heine, who developed what are now called basic hypergeometric series, the $q$-analogues of the Gauss hypergeometric series \cite{Heine1847,GasperRahman2004}. The systematic study of $q$-difference operators and $q$-integration began with F.~H.~Jackson, who defined the $q$-derivative and the $q$-definite integral in the early twentieth century \cite{Jackson1908,Jackson1910}. Soon afterwards the theory of linear $q$-difference equations was developed by Carmichael \cite{Carmichael1912}, Mason \cite{Mason1915}, Adams \cite{Adams1929} and Trjitzinsky \cite{Trjitzinsky1933}, and $q$-series also appear prominently in the work of Ramanujan. \\

For $0 < q < 1$, the $q$-difference operator (or Jackson $q$-derivative) of an analytic function $f$ is defined by
$$
\partial_q f(z) =  \begin{cases}  \frac{f(z) - f(qz)}{(1 - q)z}, & z \neq 0, \\ f'(0), & z = 0. \end{cases}
$$
Observe that $\lim_{q \to 1^-} \partial_q f(z) = f'(z)$, and $\lim_{q \to 0^+} \partial_q f(z) = \frac{f(z)}{z}$ whenever $f(0) = 0$. Since Jackson's pioneering work \cite{Jackson1908,Jackson1910}, $q$-calculus has evolved into a prominent framework across analysis and geometric function theory. From the definition of $\partial_q$, for any $f, g \in \mathcal{A}$ and $z \in \mathbb{D}$, the operator satisfies
\begin{enumerate}
\item[(i)] $\partial_q(f(z) \pm g(z)) = \partial_q(f(z)) \pm \partial_q(g(z))$,\\
\item[(ii)] $\partial_q(f(z)g(z)) = g(z)\partial_q(f(z)) + f(qz)\partial_q(g(z))$,\\
\item[(iii)] $\partial_q\left(\frac{f(z)}{g(z)}\right) = \frac{g(z)\partial_q(f(z)) - f(z)\partial_q(g(z))}{g(z)g(qz)}, \quad g(z)g(qz) \neq 0$.
\end{enumerate}
For any $f \in \mathcal{A}$ of the form \eqref{FL-001}, the $q$-difference operator $\partial_q$ of $f$ can be written as
$$\partial_q f(z) = \frac{f(z) - f(qz)}{(1-q)z} = \sum_{n=1}^{\infty} [n]_q a_n z^{n-1},$$
where
$$[n]_q = \frac{1-q^n}{1-q} = 1 + q + q^2 + \dots + q^{n-1}$$
is the standard $q$-bracket (or $q$-integer). Note that $\lim_{q \to 1^-} [n]_q = n$. Correspondingly, the $q$-integral (see \cite{Jackson1910}) of a function $f \in \mathcal{A}$ is defined by
$$\int_0^z f(t) \, d_q t = z(1-q) \sum_{k=0}^{\infty} q^k f(zq^k), \qquad z \in \mathbb{D},$$
provided the series converges. For $f(z) = \sum_{n=0}^{\infty} a_n z^n$, we have
$$
\int_0^z f(t) \, d_q t = \sum_{n=0}^{\infty} \frac{a_n}{[n+1]_q} z^{n+1}.
$$
This operator satisfies the fundamental relationship $\partial_q \left( \int_0^z f(t) \, d_q t \right) = f(z)$, and as $q \to 1^-$, it reduces to the classical integral $\int_0^z f(t) \, dt$.\\

For $0 < q < 1$, the concepts of starlikeness and convexity can be generalized using the $q$-difference operator. A function $f \in \mathcal{A}$ is called $q$-starlike if
$$\operatorname{Re} \left( \frac{z \, \partial_q f(z)}{f(z)} \right) > 0, \qquad z \in \mathbb{D},$$
and $q$-convex if
$$\operatorname{Re} \left( \frac{\partial_q (z \, \partial_q f(z))}{\partial_q f(z)} \right) > 0, \qquad z \in \mathbb{D}.$$

In geometric function theory, the first link with $q$-theory was made by Ismail et al. \cite{IMS1990}, who introduced the class of $q$-starlike functions by replacing the ordinary derivative with the $q$-difference operator. This direction has since been developed by many authors; see, for instance, \cite{AgrawalSahoo2014,AgrawalSahoo2017,Srivastava2020,PiejkoSokolTrabkaWieclaw2019,PiejkoSokol2020,SelvakumaranChoiPurohit2021,Ucar2016, VermaKumarSokol2022}, and \cite{ArlGupAgar2013} for applications of $q$-calculus in operator theory.\\

For $f, g \in \mathcal{H}$, $f$ is said to be subordinate to $g$, written $f \prec g$, if there exists a Schwarz function $\omega: \mathbb{D} \to \mathbb{D}$ (i.e., analytic with $\omega(0)=0$ and $\vert{}\omega(z)\vert{} < 1$ for $z \in \mathbb{D}$) such that $f(z) = g(\omega(z))$ for all $z \in \mathbb{D}$. Notably, if $g$ is univalent in $\mathbb{D}$, then $f \prec g$ if and only if $f(0) = g(0)$ and $f(\mathbb{D}) \subseteq g(\mathbb{D})$.\\

The general theory of differential subordinations, extensively developed by Miller and Mocanu \cite{Miller-Mocanu-book-2000}, relies centrally on boundary properties of analytic mappings. A cornerstone in this study is the classical lemma due to Jack \cite{Jack-1971}.

\begin{lem}[\textbf{Jack's Lemma} \cite{Jack-1971}]\label{jack-1}
Let $\omega \in \mathcal{H}$ with $\omega(0) = 0$ and $\omega(z) \not\equiv 0$. If $|\omega(z)|$ attains its maximum value on the closed disk $|z| \leq r$ at a point $z_0$ with $|z_0| = r$, then
\[
\frac{z_0 \omega'(z_0)}{\omega(z_0)} = k,
\]
for some real number $k \geq 1$.
\end{lem}

The classical Jack lemma  was further extended by Miller and Mocanu \cite{miller-1978} and this can be seen as a higher-order extension of Jack lemma.

\begin{lem}[Miller--Mocanu Lemma \cite{miller-1978}] \label{jack-2}
Let $z_0 \in \mathbb{D}$ and $r_0 = |z_0|$. Let
\[
\omega(z) = a_n z^n + a_{n+1} z^{n+1} + a_{n+2} z^{n+2} + \dots,\quad n\ge1,
\]
be continuous on $\overline{\mathbb{D}}_{r_0}$ and analytic on $\mathbb{D}_{r_0} \cup \{z_0\}$ with $\omega(z) \not\equiv 0$. If
\[
|\omega(z_0)| = \max \{ |\omega(z)| : z \in \overline{\mathbb{D}}_{r_0} \},
\]
then there exists a real number $m \geq n$ such that
\[
\frac{z_0 \omega'(z_0)}{\omega(z_0)} = m.
\]
\end{lem}

Jack's lemma (Lemma~\ref{jack-1}) and its generalization (Lemma~\ref{jack-2}) are
central tools in the theory of differential subordination
(see~\cite{Miller-Mocanu-book-2000}). Their counterparts in Jackson's $q$-calculus
remain largely unexplored, since the quotient and chain rules of $q$-calculus differ
from the classical ones. In this article we establish a $q$-analogue of Jack's lemma.
Unlike the classical case, where the boundary multiplier is real, its $q$-analogue
$k_q$ is in general complex, with $\operatorname{Re}k_q\ge[n]_q$ sharp. As
consequences, we obtain boundary growth estimates, a $q$-extension of Nunokawa's
theorem on the imaginary axis, a $q$-analogue of the Miller--Mocanu lemma with an
explicit disk localization and a directed secant relation, and sharp circular
enclosure bounds for first-order $q$-differential subordinations. Each result
recovers its classical analogue as $q\to1^{-}$.

\section{$q$-Jack lemma}
We extend Lemma \ref{jack-2} to the $q$-calculus setting by employing the Jackson $q$-derivative. The following result characterizes the behavior of an analytic function at a point where its modulus attains a maximum on the boundary, yielding a $q$-differential inequality with a sharp lower bound expressed in terms of a $q$-bracket.

\begin{lem}[\textbf{$q$-Jack Lemma}]\label{FL-t-001}
Let $z_0 \in \mathbb{D}$ and $r_0 = |z_0| > 0$. Let  
\[
\omega(z) = \sum_{k=n}^{\infty} a_k z^k, \quad (n \ge 1),
\]
be continuous on $\overline{\mathbb{D}}_{r_0}$ and analytic on $\mathbb{D}_{r_0} \cup \{z_0\}$ with $\omega(z) \not\equiv 0$. If  
\[
|\omega(z_0)| = \max\{|\omega(z)| : z \in \overline{\mathbb{D}}_{r_0}\},
\]
then there exists a complex number $k_q$ such that  
\[
z_0 \, \partial_q \omega(z_0) = k_q \, \omega(z_0), \quad \text{where } \operatorname{Re}(k_q) \ge [n]_q = \frac{1 - q^n}{1 - q}.
\]
Furthermore, $\operatorname{Re}(k_q) = [n]_q$ if and only if $\omega(z) = C z^n$ for some non-zero constant $C \in \mathbb{C}$.
\end{lem}

\begin{proof}
Since $\omega(z) \not\equiv 0$, the maximum modulus principle implies that $|\omega(z_0)| = \max_{|z| \le r_0} |\omega(z)| > 0$. Consider the function
\[
\psi(z) = \frac{\omega(z)}{z^n}.
\]
Because $\omega(z)$ has a zero of order at least $n~ (\ge 1)$ at the origin, the singularity of $\psi$ at $z=0$ is removable. Hence, $\psi$ is continuous on $\overline{\mathbb{D}}_{r_0}$ and analytic on $\mathbb{D}_{r_0} \cup \{z_0\}$. On the boundary circle $|z| = r_0$, we have
\[
|\psi(z)| = \frac{|\omega(z)|}{r_0^n} \le \frac{|\omega(z_0)|}{r_0^n} = |\psi(z_0)|, \quad (|z| = r_0).
\]
By the maximum modulus principle applied to $\psi$ on $\overline{\mathbb{D}}_{r_0}$,
\[
|\psi(z)| \leq |\psi(z_0)| \quad \text{for all } |z| \leq r_0.
\]
In particular, evaluating at $z = q z_0$ for $0 < q < 1$ (which satisfies $|q z_0| = q r_0 < r_0$) yields
\[
\left| \frac{\psi(q z_0)}{\psi(z_0)} \right| \leq 1,
\]
which implies that
\begin{equation}\label{FL-010}
\left| \frac{\omega(q z_0)}{\omega(z_0)} \right| = \left| q^n \frac{\psi(q z_0)}{\psi(z_0)} \right| \leq q^n.
\end{equation}
From the definition of the $q$-difference operator $\partial_q$, we have
\[
k_q := \frac{z_0 \, \partial_q \omega(z_0)}{\omega(z_0)} = \frac{1}{1 - q} \left(1 - \frac{\omega(q z_0)}{\omega(z_0)}\right).
\]
Taking the real part and applying \eqref{FL-010}, we obtain
\[
\operatorname{Re}(k_q) = \frac{1}{1 - q}\left(1 - \operatorname{Re}\left(\frac{\omega(q z_0)}{\omega(z_0)}\right)\right) \ge \frac{1}{1 - q}\left(1 - \left|\frac{\omega(q z_0)}{\omega(z_0)}\right|\right) \ge \frac{1 - q^n}{1 - q} = [n]_q.
\]

To analyze the case of equality, suppose $\operatorname{Re}(k_q) = [n]_q$. This equality requires both 
$$\operatorname{Re}\left(\frac{\omega(q z_0)}{\omega(z_0)}\right) = \left|\frac{\omega(q z_0)}{\omega(z_0)}\right|\quad\text{and}~~ \left|\frac{\omega(q z_0)}{\omega(z_0)}\right| = q^n.$$ The second condition implies that $|\psi(q z_0)| = |\psi(z_0)|$. Since $q z_0 \in \mathbb{D}_{r_0}$ is an interior point and $|\psi|$ attains its maximum on $\overline{\mathbb{D}}_{r_0}$ at $q z_0$, the strong maximum principle asserts that $\psi(z)$ must be identically constant on $\overline{\mathbb{D}}_{r_0}$. Thus, $\omega(z) = C z^n$ for some constant $C \neq 0$.

Conversely, if $\omega(z) = C z^n$ with $C \neq 0$, direct computation yields $\frac{\omega(q z_0)}{\omega(z_0)} = q^n$, whence $k_q = \frac{1 - q^n}{1 - q} = [n]_q$, and therefore $\operatorname{Re}(k_q) = [n]_q$.
\end{proof}

%%%%%%%%%%%%%%%%%%%%%%%%%%%%%%%%%%%%%%%%%
\begin{rem}\label{FL-tr-001}
In view of Lemma \ref{FL-t-001} ($q$-Jack's Lemma), we have
\[
k_q = \frac{z_0 \, \partial_q \omega(z_0)}{\omega(z_0)} = \frac{\omega(z_0) - \omega(q z_0)}{(1-q)\,\omega(z_0)}.
\]
Since $\omega$ is analytic at $z_0$, its Taylor expansion about $z_0$ along the radial segment gives
\[
\omega(q z_0) = \omega(z_0) + \omega'(z_0)(q-1)z_0 + \mathcal{O}((1-q)^2) \quad \text{as } q \to 1^-.
\]
Substituting this expansion into the expression for $k_q$ yields
\[
k_q = \frac{z_0 \omega'(z_0)}{\omega(z_0)} + \mathcal{O}(1-q).
\]
Taking the limit as $q \to 1^-$ and recalling that $\lim_{q \to 1^-} [n]_q = n$, we deduce that
\[
k := \lim_{q \to 1^-} k_q = \frac{z_0 \omega'(z_0)}{\omega(z_0)}.
\]
Moreover, because $\operatorname{Re}(k_q) \ge [n]_q$ for each $q \in (0,1)$, it follows that $\operatorname{Re}(k) \ge n$. Since the classical Miller--Mocanu lemma (Lemma \ref{jack-2}) ensures that $\frac{z_0 \omega'(z_0)}{\omega(z_0)}$ is real, we conclude that $\lim_{q \to 1^-} \operatorname{Im}(k_q) = 0$ and $k \ge n$, thereby recovering the classical result.
\end{rem}

A remarkable qualitative distinction between the classical Jack lemma and its $q$-analogue is that while the classical parameter $k = \frac{z_0 \omega'(z_0)}{\omega(z_0)}$ is strictly real, the quantity $k_q$ in Lemma \ref{FL-t-001} can be genuinely non-real for every $q \in (0, 1)$. The following explicit construction justifies this phenomenon at a boundary point lying strictly within the open unit disk $\mathbb{D}$.

\begin{example}\label{FL-tr-005}
Let $z_0 = \frac{1}{2} \in \mathbb{D}$ on the circle $|z| = r_0 = \frac{1}{2}$, so that the closed disk $\overline{\mathbb{D}}_{r_0} \subset \mathbb{D}$. Consider the polynomial $\omega$ of degree $3$ (with order of vanishing $n = 1$ at $z = 0$) defined by
\begin{equation*}
\omega(z) = \left(\frac{1}{2} - \frac{i}{5}\right) z + \left(2 + \frac{4i}{5}\right) z^2 + \left(2 - \frac{4i}{5}\right) z^3.
\end{equation*}
A direct computation shows that
\[
\omega(z_0) = \omega\left(\frac{1}{2}\right) = \left(\frac{1}{4} + \frac{1}{2} + \frac{1}{4}\right) + i\left(-\frac{1}{10} + \frac{2}{10} - \frac{1}{10}\right) = 1,
\]
and so $|\omega(z_0)| = 1$. 
To verify that $|\omega(z)| \le 1$ for all $z \in \overline{\mathbb{D}}_{1/2}$, write $z = \frac{1}{2} e^{i\theta}$ with $\theta \in [-\pi, \pi]$ on the boundary circle $|z| = \frac{1}{2}$. Thus
\[
\begin{aligned}
\omega\left(\frac{1}{2} e^{i\theta}\right) &= e^{2i\theta} \left[ \left(\frac{1}{4} - \frac{i}{10}\right)e^{-i\theta} + \left(\frac{1}{2} + \frac{i}{5}\right) + \left(\frac{1}{4} - \frac{i}{10}\right)e^{i\theta} \right] \\
&= e^{2i\theta} \left[ \left(\frac{1}{2} + \frac{i}{5}\right) + 2\left(\frac{1}{4} - \frac{i}{10}\right)\cos\theta \right] \\
&= e^{2i\theta} \left[ \frac{1 + \cos\theta}{2} + i\,\frac{1 - \cos\theta}{5} \right].
\end{aligned}
\]
Taking the squared modulus yields
\[
\left|\omega\left(\frac{1}{2} e^{i\theta}\right)\right|^2 = \frac{(1+\cos\theta)^2}{4} + \frac{(1-\cos\theta)^2}{25} =\frac{(1+u)^2}{4} + \frac{(1-u)^2}{25}=:g(u),
\]
where $u = \cos\theta \in [-1, 1]$.
Since $g''(u) = \frac{29}{50} > 0$, the function $g$ is strictly convex on $[-1, 1]$, so its global maximum on $[-1, 1]$ occurs at an endpoint:
\[
g(1) = 1 \quad (\theta = 0), \qquad g(-1) = \frac{4}{25} = 0.16 \quad (\theta = \pm \pi).
\]
By the maximum modulus principle, it follows that
\[
\max_{|z| \le 1/2} |\omega(z)| =1= |\omega(z_0)|.
\]

Now, for any $q \in (0, 1)$, evaluating $\omega$ at $q z_0 = \frac{q}{2}$ gives
\[
\begin{aligned}
\omega(q z_0) = \omega\left(\frac{q}{2}\right) &= \left(\frac{1}{4} - \frac{i}{10}\right) q + \left(\frac{1}{2} + \frac{i}{5}\right) q^2 + \left(\frac{1}{4} - \frac{i}{10}\right) q^3 \\
&= \frac{q(1+q)^2}{4} - \frac{i}{10} q(1-q)^2.
\end{aligned}
\]
Since $\omega(z_0) = 1$, the quantity $k_q$ is given by
\[
k_q = \frac{z_0 \, \partial_q \omega(z_0)}{\omega(z_0)} = \frac{1}{1 - q}\left(1 - \frac{\omega(q z_0)}{\omega(z_0)}\right) = \frac{1 - \omega\left(\frac{q}{2}\right)}{1 - q}.
\]
Separating into real and imaginary components, we have
\[
\operatorname{Re}(k_q) = \frac{1 - \frac{q(1+q)^2}{4}}{1 - q} = 1 + q + \frac{q^2}{4} > [1]_q = 1,
\]
and
\[
\operatorname{Im}(k_q) = -\frac{\operatorname{Im}\left(\omega\left(\frac{q}{2}\right)\right)}{1 - q} = \frac{1}{10} q (1 - q).
\]
Thus, for every $q \in (0, 1)$, $\operatorname{Im}(k_q) > 0$, confirming that $k_q \notin \mathbb{R}$. As expected from Remark \ref{FL-tr-001}, $\lim_{q \to 1^-} \operatorname{Im}(k_q) = 0$, recovering the real limit $k = \lim_{q \to 1^-} k_q = \frac{9}{4} \ge 1$.
\end{example}

%%%%%%%%%%%%%%%%%%%%%%%%%%%%%%%%%%%%%%%%%%%%%%%%%%%

\section{Applications of $q$-Jack lemma}
The classical Jack lemma \cite{Jack-1971} and its extension by Miller and Mocanu \cite{miller-1978} serve as the foundational technical engine in geometric function theory for establishing differential subordinations, geometric inclusion relations, and coefficient bounds. In this section, we demonstrate the utility of the $q$-Jack lemma (Lemma~\ref{FL-t-001}) across three primary geometric contexts:
\begin{enumerate}[(i)]
\item establishing boundary growth and mapping behavior for functions not bounded by the unit disk;
\item generalizing Nunokawa's boundary theorem on the vanishing of the real part;
\item developing $q$-analogues of the Miller--Mocanu boundary subordination lemma and first-order differential enclosures.
\end{enumerate}

We begin by applying Lemma~\ref{FL-t-001} to characterize the behavior of analytic functions mapping outside the closed unit disk.

\begin{thm}\label{FL-t-005}
Let $f$ be a non-constant analytic function in the unit disk $\mathbb{D}$ with $f(z) = a_n z^n + a_{n+1} z^{n+1} + \dots$ ($n \ge 1$, $a_n \neq 0$) and $f(0) = 0$. If $f(\mathbb{D}) \not\subseteq \overline{\mathbb{D}}$ (that is, $f(\mathbb{D})$ contains points outside the closed unit disk $\overline{\mathbb{D}}$), then
\begin{enumerate}[(i)]
\item there exists a point $z_0 \in \mathbb{D}$ such that $|\partial_q f(z_0)| > [n]_q$;
\item there exists a point $z_1 \in \mathbb{D}$ and a complex constant $k_q$ such that $\vert{}f(z_1)\vert{} = 1$ and
$$\frac{z_1 \, \partial_q f(z_1)}{f(z_1)} = k_q, \quad \text{where } \operatorname{Re}(k_q) \ge [n]_q = \frac{1 - q^n}{1 - q}.
$$
\end{enumerate}
\end{thm}

\begin{proof}
Define $M(\rho) := \max_{\vert{}z\vert{}=\rho} \vert{}f(z)\vert{}$ for $0 < \rho < 1$, and $M := \sup_{z \in \mathbb{D}} \vert{}f(z)\vert{}$.
Since $f(\mathbb{D})$ contains points outside $\overline{\mathbb{D}}$, we have $M > 1$. By the maximum modulus principle, $M(\rho)$ is strictly increasing and continuous on $(0,1)$, so there exists $r_0 \in (0,1)$ such that
$$
M(\rho) > 1 \quad \text{for all } \rho \in (r_0, 1).
$$
For each $r \in (r_0, 1)$, choose $z_r \in \partial \mathbb{D}_r = \{z \in \mathbb{C} : \vert{}z\vert{} = r\}$ such that $\vert{}f(z_r)\vert{} = M(r)$.
Applying Lemma \ref{FL-t-001} ($q$-Jack's Lemma) to $f(z)$ at $z_r$, there exists $k_q \in \mathbb{C}$ such that
$$
\frac{z_r \, \mathbf{\partial_q} f(z_r)}{f(z_r)} = k_q, \quad \text{where } \operatorname{Re}(k_q) \ge [n]_q.
$$
Taking the absolute value yields
$$\left\vert{} \frac{z_r \, \partial_q f(z_r)}{f(z_r)} \right\vert{} = \vert{}k_q\vert{} \ge \operatorname{Re}(k_q) \ge [n]_q \ge 1.
$$
Therefore,
$$\vert{}\partial_q f(z_r)\vert{} \ge \frac{\vert{}f(z_r)\vert{}}{\vert{}z_r\vert{}} [n]_q = \frac{M(r)}{r} [n]_q > \frac{1}{r} [n]_q > [n]_q.
$$
Setting $z_0 = z_r$ completes the proof of part (i).\\

To prove part (ii), observe that $M(0) = 0$ and $\lim_{r \to 1^-} M(r) = M > 1$. Because $M(r)$ is continuous and strictly increasing, the Intermediate Value Theorem ensures the existence of a unique radius $r_1 \in (0,1)$ such that $M(r_1) = 1$. By the monotonicity of $M(r)$, we have
\[
|f(z)| \le M(r_1) = 1 \quad \text{for all } |z| \le r_1.
\]
Choosing $z_1 \in \partial \mathbb{D}_{r_1}$ such that $|f(z_1)| = M(r_1) = 1$, the function $|f|$ attains its maximum on the closed disk $\overline{\mathbb{D}}_{r_1}$ at the boundary point $z_1$. Applying Lemma \ref{FL-t-001} to $f$ on $\overline{\mathbb{D}}_{r_1}$ at $z_1$, there exists a complex number $k_q$ such that
\[
\frac{z_1 \, \partial_q f(z_1)}{f(z_1)} = k_q, \quad \text{where } \operatorname{Re}(k_q) \ge [n]_q,
\]
which establishes part (ii).

\end{proof}

%%%%%%%%%%%%%%%%%%%%%%%%%%%%%%%%%%%%%%%%%%
Another fundamental cornerstone in geometric function theory is the analysis of boundary points at which the real part of an analytic function first vanishes. In the classical setting, this behavior is governed by the following well-known result due to Nunokawa \cite{nunokawa-1992}.

\begin{customthm}{A}\cite{nunokawa-1992}\label{FL-t-010}
Let $p(z)$ be analytic in $\mathbb{D}$ with $p(0) = 1$. Suppose that there exists a point $z_{0} \in \mathbb{D}$ such that $\operatorname{Re} p(z_0) = 0$, $p(z_0) \neq 0$ and $\operatorname{Re} p(z) > 0$ for $\vert{}z\vert{} < \vert{}z_0\vert{}.$ Then at $z = z_0$, we have
\begin{equation}\label{FL-020}
\frac{z_0 p'(z_0)}{p(z_0)} = i k,
\end{equation}
where $k$ is a real number satisfying $\vert{}k\vert{} \ge 1$.\end{customthm}

Theorem \ref{FL-t-010} plays a prominent role in geometric function theory, particularly in deriving criteria for starlikeness, convexity, and strongly starlike functions. Geometrically, condition \eqref{FL-020} reflects the fact that as $z$ approaches $z_0$ along the circle $|z| = |z_0|$, the image curve $p(|z_0|e^{i\theta})$ makes contact with the imaginary axis, forcing the logarithmic derivative to point strictly along the imaginary direction.\\

By employing the $q$-Jack lemma established in Lemma \ref{FL-t-001}, we now derive a basic $q$-analogue of Theorem \ref{FL-t-010}. This extension describes the action of the Jackson $q$-difference operator at the boundary of the region of positivity.

\begin{thm}\label{FL-t-015}
Let $p(z)$ be analytic in $\mathbb{D}$ with $p(0) = 1$. Suppose there exists a point $z_0 \in \mathbb{D}$ such that $\operatorname{Re} p(z_0) = 0$, $p(z_0) \neq 0$ and $\operatorname{Re} p(z) > 0$ for $|z| < |z_0|$. Then at $z = z_0$, there exists a complex constant $\lambda_q$ such that
$$
\frac{z_0 \, \partial_q p(z_0)}{p(z_0)} = \mathbf{-} \frac{\lambda_q (1 - p(z_0))(1 + p(qz_0))}{2 p(z_0)}, \quad \text{where } \operatorname{Re}(\lambda_q) \ge 1.
$$
\end{thm}

\begin{proof}
Consider the M\"obius transformation
\begin{align}\label{FL-025}
\phi(z) = \frac{1 - p(z)}{1 + p(z)}.
\end{align}
Then $\phi(0) = 0$  and $|\phi(z)| < 1$ for $|z| < |z_0|$. Since $\operatorname{Re} p(z_0) = 0$ and $p(z_0) \neq 0$, we can write $p(z_0) = i a$ with $a \in \mathbb{R} \setminus \{0\}$, and so
\[
|\phi(z_0)| = \left| \frac{1 - i a}{1 + i a} \right| = 1.
\]
Hence $|\phi|$ attains its maximum on the circle $|z| = |z_0|$ at $z_0$. Applying Lemma \ref{FL-t-001} ($q$-Jack's Lemma) to $\phi(z)$ with $n = 1$, there exists a complex number $\lambda_q$ such that
\begin{align}\label{FL-030}
    z_0 \parq \phi(z_0) = \lambda_q \phi(z_0), \quad \operatorname{Re}\lambda_q \ge 1.
\end{align}
From \eqref{FL-025}, using the $q$-quotient rule, we get
$$
\partial_q \phi(z) = \frac{(1 + p(z))(-\partial_q p(z)) - (1 - p(z))\partial_q p(z)}{(1 + p(z))(1 + p(qz))} = \mathbf{-} \frac{2 \partial_q p(z)}{(1 + p(z))(1 + p(qz))}.
$$
Evaluating at $z = z_0$ and substituting into \eqref{FL-030}, we get
$$
- \frac{2 z_0 \, \partial_q p(z_0)}{(1 + p(z_0))(1 + p(qz_0))} = \lambda_q \cdot \frac{1 - p(z_0)}{1 + p(z_0)}
$$
and so
\begin{align}\label{FL-035}
\frac{z_0 \parq p(z_0)}{p(z_0)} =- \frac{\lambda_q (1 - p(z_0))(1 + p(qz_0))}{2 p(z_0)},
\end{align}
which completes the proof.
\end{proof}

\begin{rem}
As $q \to 1^-$, we have $\partial_q p(z_0) \to p'(z_0)$ and $p(q z_0) \to p(z_0) = i a$. Moreover, by the classical Jack lemma, the complex parameter $\lambda_q$ satisfies
\[
\lim_{q \to 1^-} \operatorname{Im}(\lambda_q) = 0 \quad \text{and} \quad \lambda := \lim_{q \to 1^-} \lambda_q \ge 1, \quad \text{with } \lambda \in \mathbb{R}.
\]
Taking the limit as $q \to 1^-$ on both sides of \eqref{FL-035}, we obtain
\[
\frac{z_0 p'(z_0)}{p(z_0)} = - \frac{\lambda (1 - i a)(1 + i a)}{2 i a} = - \frac{\lambda (1 + a^2)}{2 i a} = i \lambda \left(\frac{1 + a^2}{2a}\right).
\]
Setting
\[
k := \lambda \left(\frac{1 + a^2}{2a}\right) \in \mathbb{R},
\]
we arrive at $\frac{z_0 p'(z_0)}{p(z_0)} = i k$. Since $|a + a^{-1}| \ge 2$ (or equivalently, $\frac{1 + a^2}{2|a|} \ge 1$) for all $a \in \mathbb{R} \setminus \{0\}$, it follows that
\[
|k| = \lambda \left(\frac{1 + a^2}{2|a|}\right) \ge 1,
\]
which directly recovers Nunokawa's Theorem \ref{FL-t-010}.
\end{rem}

%%%%%%%%%%%%%%%%%%%%%%%%%%%%%%%%%%%%%%%%%
Next, we turn our attention to general subordination relationships. In the systematic framework of Miller and Mocanu \cite{Miller-Mocanu-book-2000}, boundary subordination theorems require target mappings to belong to an admissible class possessing boundary regularity. To this end, we denote by $\mathcal{Q}$ the class of functions $q$ that are analytic and injective on $\overline{\mathbb{D}} \setminus E(q)$, where
\[
E(q) = \left\{ \zeta \in \partial \mathbb{D} : \lim_{z \to \zeta} q(z) = \infty \right\},
\]
and satisfy $q'(\zeta) \neq 0$ for all $\zeta \in \partial \mathbb{D} \setminus E(q)$. In differential subordination theory, this class serves as the standard collection of univalent target functions. The following lemma provides a fundamental $q$-analogue of the celebrated Miller--Mocanu boundary subordination lemma \cite[Lemma 2.2d, p. 24]{Miller-Mocanu-book-2000}.

%%%%%%%%%%%%%%%%%%%%%%%%%%%%%%%%
\begin{thm}\label{FL-t-020}
Let $g \in \mathcal{Q}$ with $g(0) = a$, and let 
\[
f(z) = a + \sum_{k=n}^{\infty} a_k z^k, \quad (n \ge 1, \; a_n \neq 0),
\]
be analytic in $\mathbb{D}$ with $f(z) \not\equiv a$. If $f$ is not subordinate to $g$ ($f \not\prec g$), then there exist points $z_0 \in \mathbb{D}$, $\zeta_0 \in \partial \mathbb{D} \setminus E(g)$, and an interior point $\zeta_1 \in \mathbb{D}$ with $|\zeta_1| \le q^n$ such that
\begin{enumerate}[(i)]
\item $f(z_0) = g(\zeta_0)$;
\item the Jackson $q$-derivative satisfies the relation
\[
z_0 \, \partial_q f(z_0) = \frac{g(\zeta_0) - g(\zeta_1)}{1 - q} = m_q \, \zeta_0 \left( \frac{g(\zeta_0) - g(\zeta_1)}{\zeta_0 - \zeta_1} \right),
\]
where $m_q = \dfrac{\zeta_0 - \zeta_1}{(1 - q)\zeta_0}$ satisfies the disk localization
\begin{equation}\label{FL-040}
\left| m_q - \frac{1}{1-q} \right| \le \frac{q^n}{1-q}.
\end{equation}
In particular, $m_q$ satisfies
\[
[n]_q = \frac{1 - q^n}{1 - q} \le \operatorname{Re}(m_q) \le \frac{1 + q^n}{1 - q}, \qquad |\operatorname{Im}(m_q)| \le \frac{q^n}{1-q},
\]
and
\[
|\arg(m_q)| \le \arcsin(q^n) < \frac{\pi}{2}.
\]
\end{enumerate}
\end{thm}

\begin{proof}
Since $f(0) = g(0) = a$ and both $f$ and $g$ are analytic in $\mathbb{D}$ with $g$ univalent, define the maximal radius
\[
r_0 := \sup \{ r \in (0, 1) : f(\mathbb{D}_r) \subset g(\mathbb{D}) \}.
\]
Because $f \not\prec g$, we have $f(\mathbb{D}) \not\subset g(\mathbb{D})$, which guarantees that $0 < r_0 < 1$. By the maximality of $r_0$, we have $f(\mathbb{D}_{r_0}) \subset g(\mathbb{D})$, $f(\overline{\mathbb{D}}_{r_0}) \subset \overline{g(\mathbb{D})}$ with the boundary contact property $f(\partial \mathbb{D}_{r_0}) \cap \partial g(\mathbb{D}) \neq \emptyset$. Therefore, there exists a point $z_0 \in \partial \mathbb{D}_{r_0}$ such that $f(z_0) \in \partial g(\mathbb{D})$.  The injectivity of $g$ on $\overline{\mathbb{D}} \setminus E(g)$ guarantees the existence of a unique $\zeta_0 \in \partial \mathbb{D} \setminus E(g)$ such that $f(z_0) = g(\zeta_0)$, proving (i).\\

Next, since $g$ is univalent on $\overline{\mathbb{D}} \setminus E(g)$, the function $\phi: \overline{\mathbb{D}}_{r_0} \to \overline{\mathbb{D}}$ given by $\phi(z) = g^{-1}(f(z))$ is well-defined and continuous in $\overline{\mathbb{D}}_{r_0}$ and satisfies $\phi(0) = g^{-1}(a) = 0$. Further, $\phi(z)$ is analytic in $\mathbb{D}_{r_0}$ and has a zero of order $n$ at the origin. 
Furthermore, $|\phi(z)| < 1$ on $\mathbb{D}_{r_0}$ and
\[
|\phi(z_0)| = |\zeta_0| = 1 = \max_{|z| \le r_0} |\phi(z)|.
\]

Setting $\zeta_1 := \phi(q z_0)$, inequality \eqref{FL-010} from Lemma \ref{FL-t-001} ensures that
\[
|\zeta_1| = |\phi(q z_0)| \le q^n |\phi(z_0)| = q^n < 1.
\]
By definition of the Jackson $q$-difference operator,
\[
m_q := \frac{z_0 \, \partial_q \phi(z_0)}{\phi(z_0)} = \frac{\phi(z_0) - \phi(q z_0)}{(1-q)\phi(z_0)} = \frac{\zeta_0 - \zeta_1}{(1-q)\zeta_0} = \frac{1}{1-q}\left(1 - \frac{\zeta_1}{\zeta_0}\right).
\]
Since $|\zeta_0| = 1$ and $|\zeta_1| \le q^n$, we immediately obtain
\[
\left| m_q - \frac{1}{1-q} \right| = \frac{1}{1-q} \left| \frac{\zeta_1}{\zeta_0} \right| \le \frac{q^n}{1-q}.
\]
This proves \eqref{FL-040}. The bounds on $\operatorname{Re}(m_q)$ and $|\operatorname{Im}(m_q)|$ follow directly from the geometry of the disk. Moreover, because the disk is centered at $\frac{1}{1-q} > 0$ with radius $\frac{q^n}{1-q} < \frac{1}{1-q}$, the maximum argument of points in this disk is given by 
$$\arcsin\left(\frac{q^n/(1-q)}{1/(1-q)}\right) = \arcsin(q^n) < \frac{\pi}{2}.$$

Finally, using $f(z) = g(\phi(z))$, we have
\[
z_0 \, \partial_q f(z_0) = \frac{g(\phi(z_0)) - g(\phi(q z_0))}{1 - q} = \frac{g(\zeta_0) - g(\zeta_1)}{1 - q} = m_q \, \zeta_0 \left( \frac{g(\zeta_0) - g(\zeta_1)}{\zeta_0 - \zeta_1} \right),
\]
which completes the proof.
\end{proof}

%%%%%%%%%%%%%%%%%%%%%%%%%%%%%%%%

\begin{rem}
As $q \to 1^-$, we have $\zeta_1 = \phi(q z_0) \to \phi(z_0) = \zeta_0$. Hence, the difference quotient satisfies
\[
\lim_{q \to 1^-} \frac{g(\zeta_0) - g(\zeta_1)}{\zeta_0 - \zeta_1} = g'(\zeta_0).
\]
Furthermore, by Remark \ref{FL-tr-001}, $\lim_{q \to 1^-} \operatorname{Im}(m_q) = 0$ and $m := \lim_{q \to 1^-} m_q \ge n$, with $m \in \mathbb{R}$. Passing to the limit in part (ii) yields
\[
z_0 f'(z_0) = m \, \zeta_0 \, g'(\zeta_0), \quad (m \ge n),
\]
which precisely recovers the classical Miller--Mocanu lemma \cite[Lemma 2.2d, p. 24]{Miller-Mocanu-book-2000}.
\end{rem}

%%%%%%%%%%%%%%%%%%%%%%%%%%%%%%
Equipped with Lemma \ref{FL-t-001} and Theorem \ref{FL-t-020}, we are now positioned to examine specific differential subordination problems.
In classical complex analysis, a widely investigated problem (see \cite{miller-1981, Miller-Mocanu-book-2000}) concerns finding the best dominant for first-order differential expressions of the form
\[
p(z) + \gamma \, z p'(z) \prec h(z).
\]
When the target domain is a disk centered at $1$, a fundamental theorem states that if $|p(z) + \gamma z p'(z) - 1| < M$, then $|p(z) - 1| < \frac{M}{1 + n\gamma}$ for functions vanishing to order $n$ at the origin. 
Our objective here is to establish the $q$-difference analogue of this circular enclosure theorem. By replacing the classical derivative with the Jackson $q$-derivative $\partial_q$, we determine how the discrete deformation parameter $q \in (0, 1)$ influences the radius of subordination, providing a natural bridge to the classical continuous result as $q \to 1^-$.

\begin{thm}\label{FL-t-025}
Let $p$ be an analytic function in the unit disk $\mathbb{D}$ of the form
\[
p(z) = 1 + \sum_{k=n}^{\infty} c_k z^k, \quad (n \ge 1).
\]
Let $\gamma \ge 0$ and $M > 0$. If
\begin{equation}\label{FL-050}
\left| p(z) + \gamma \, z \, \partial_q p(z) - 1 \right| < M, \quad (z \in \mathbb{D}),
\end{equation}
then
\[
|p(z) - 1| < \frac{M}{1 + \gamma [n]_q}, \quad (z \in \mathbb{D}).
\]
\end{thm}

\begin{proof}
Consider the function $\omega: \mathbb{D} \to \mathbb{C}$ defined by $\omega(z) := p(z) - 1.$
Then $\omega(0) = 0$ and $\omega$ has a zero of order at least $n$ at the origin.
To establish the theorem, it suffices to prove that $|\omega(z)| < R$ for all $z \in \mathbb{D}$, where
\[
R := \frac{M}{1 + \gamma [n]_q} > 0.
\]

Suppose, to the contrary, that there exists an interior point in $\mathbb{D}$ where $|\omega(z)| \ge R$. Since $\omega(0) = 0 < R$ and $\omega$ is analytic in $\mathbb{D}$, the maximum modulus principle ensures the existence of a radius $r_0 \in (0, 1)$ and a boundary point $z_0 \in \partial \mathbb{D}_{r_0}$ such that
\[
|\omega(z_0)| = \max_{|z| \le r_0} |\omega(z)| = R.
\]
Applying Lemma \ref{FL-t-001} ($q$-Jack's Lemma) to $\omega$ on the closed disk $\overline{\mathbb{D}}_{r_0}$, there exists a complex constant $k_q \in \mathbb{C}$ such that
\begin{equation}\label{FL-060}
z_0 \, \partial_q \omega(z_0) = k_q \, \omega(z_0), \quad \text{where } \operatorname{Re}(k_q) \ge [n]_q.
\end{equation}
Evaluating the differential expression in \eqref{FL-050} at $z = z_0$ and using \eqref{FL-060}, we obtain
\[
\begin{aligned}
\left| p(z_0) + \gamma \, z_0 \, \partial_q p(z_0) - 1 \right| &= \left| \omega(z_0) + \gamma \, z_0 \, \partial_q \omega(z_0) \right| \\
&= \left| \omega(z_0) + \gamma \, k_q \, \omega(z_0) \right| \\
&= |\omega(z_0)| \cdot |1 + \gamma k_q| \\
&= R \, |1 + \gamma k_q|\\
&\ge R(1 + \gamma \operatorname{Re}(k_q)) \\
&\ge R \left( 1 + \gamma[n]_q \right)=M,
\end{aligned}
\]
which contradicts the hypothesis \eqref{FL-050}. Consequently, no such point $z_0$ can exist, and we conclude that
\[
|p(z) - 1| < \frac{M}{1 + \gamma[n]_q} \quad \text{for all } z \in \mathbb{D}.
\]
\end{proof}

\noindent\textbf{Data availability:}
Data sharing not applicable to this article as no data sets were generated or analyzed during the current study.\vspace{4mm}
\\
\noindent\textbf{Authors Contributions:}
All authors contributed equally to the investigation of the problem and the order of the authors is given alphabetically according to their surname. All authors read and approved the final manuscript.\vspace{4mm}
\\
\noindent\textbf{Acknowledgement:}
The second author thanks the CSIR for the financial support through CSIR-SRF Fellowship ( File no.: 09/0973(13731)/2022-EMR-I ).

\end{document}